\documentclass [11pt]{article}
\usepackage{amssymb,amsmath,comment,amsthm}

\def\N{\mathbb{N}}

\def\F{\mathbb{F}}

\newtheorem{theorem}{Theorem}[section]

\newtheorem{corollary}[theorem]{Corollary}
\newtheorem{lemma}[theorem]{Lemma}
\newtheorem{definition}[theorem]{Definition}
\newtheorem{remark}[theorem]{Remark}

\begin{document}
\title{On special perfect polynomials  over $\F_2$}
 \author{Olivier Rahavandrainy\\
Univ Brest, UMR CNRS 6205\\
Laboratoire de Math\'ematiques de Bretagne Atlantique}
\maketitle
\\
{\bf{Abstract}}
We give some necessary conditions satisfied by one family of special perfect polynomials
over $\F_2$.
{\section{Introduction}}
In this paper, we work over the ground field $\F_2 = \{0,1\}$ and over its first extension $\F_4 = \{0,1,\alpha, \alpha + 1\}$, where
$\alpha^2+\alpha +1 = 0$.\\
As usual, $\N$ (resp. $\N\sp{*}$) denotes the set of nonnegative
integers (resp. of positive integers).

We say that $A \in \F_2[x]$ is odd (resp. even) if $\gcd(A,x^2+x) = 1$ (resp. $\not= 1$).
We denote by $\sigma(A)$ (resp. $\omega(A)$)  the sum of all divisors  (resp. the number of all irreducible divisors) of $A$. The function $\sigma$ is multiplicative.
We say that $A$ is \emph{perfect} if $\sigma(A) = A$.

We write, for $m \in \N^*$ and for an irreducible divisor $S$ of $A$:
$$\mbox{$S^m \| A$ if $S^m$ divides $A$ but $S^{m+1}$ does not divide it.}$$
Finally, we say (\cite{Gall-Rahav-up-Mersenne}) that a binary polynomial is a Mersenne one if it is of the form $x^a(x+1)^b +1$ for some positive coprime integers $a,b$.

It is easy to see that every odd perfect polynomial (over $\F_2$) is a square. We say that $A$ is ``{\it{special}}'' (\cite{Gall-Rahav-Bdx}) if it is odd and if each exponent (involving in $A$) is equal to $2$: $A =\displaystyle{\prod_i {T_i}^2}$, where each $T_i$ is irreducible.

Several studies are done about perfect polynomials over $\F_2$ namely, the existence or not of the odd ones (see \cite{Canaday}, \cite{Gall-Rahav-F4}, \cite{Gall-Rahav-Bdx}, \cite{Gall-Rahav8}, \cite{Gall-Rahav-up-Mersenne} and references therein). In particular (\cite{Gall-Rahav-Bdx}, Theorem 5.5), if $A$ is special then
$\omega(A) \geq 10$ and $\deg(S) \geq 30$, for any irreducible factor $S$ of $A$.

If $A$ is special and perfect, then we may write:
\begin{equation} \label{forme-de-A}
A = {P_1}^2 \cdots P_{n_p}^{2} \ {Q_1}^2 \cdots Q_{n_q}^2 \ {R_1}^2 \cdots R_{n_r}^2,\\
\mbox{where $n_p, n_q, n_r \in \N^*$.}
\end{equation}
By Lemma \ref{mini-medium-maxi}, one has for $1\leq i,j\leq n_p$, $1 \leq k, \ell \leq n_r$:
$$\deg(P_i) = \deg(P_{j})  < \deg(Q_1) \leq \cdots \leq \deg(Q_{n_q}) < \deg(R_k) = \deg(R_{\ell}).$$
We call (Definition \ref{min-med-max}) the $P_i$'s: minimal (resp. the $Q_j$'s:  medium  and the $R_k$'s:  maximal)  prime divisors of $A$.

We consider the sets of minimal, medium and maximal prime divisors of $A$ and their degrees:
\begin{equation} \label{lesPQR}
\begin{array}{l}
{\cal{M}}_{in}: =  \{P_1,\ldots, P_{n_p}\}, {\cal{M}}_{ed}: =  \{Q_1,\ldots, Q_{n_q}\}, {\cal{M}}_{ax}: =  \{R_1,\ldots, R_{n_r}\}\\
p:= \deg(P), \mbox{ for $P \in {\cal{M}}_{in}$}, \ q_j:=\deg(Q_j) \mbox{ and $r:= \deg(R),$ for $R \in {\cal{M}}_{ax}$}.
\end{array}
\end{equation}
Moreover, the degrees $p, r$ and $q_1,\ldots,q_r$ are all even (Lemma \ref{alldegre-even}).\\

Set $q = \displaystyle{\max_j{q_j}}$,
${\cal{D}}_A:= {\cal{M}}_{in} \cup {\cal{M}}_{ed} \cup {\cal{M}}_{ax}$ and
\begin{equation} \label{lesPQR2}
\begin{array}{l}
{\cal{M}}_{{in}1}: = \{P \in {\cal{M}}_{in}: \sigma(P^2) \mbox{ is irreducible}\},\
{\cal{M}}_{{in}2}: = {\cal{M}}_{in} \setminus {\cal{M}}_{{in}1},\\
{\cal{S}}_{1}: = \{\sigma(P^2): P \in {\cal{M}}_{{in}1}\}.
%{\cal{M}}_{{ed}1}: = \{\sigma(P^2): P \in {\cal{M}}_{{in}1}\},\
%{\cal{M}}_{{ed}2}: = {\cal{M}}_{ed} \setminus {\cal{M}}_{{ed}1}.
\end{array}
\end{equation}
In this paper, we give several necessary conditions satisfied by a family of special perfect polynomials (if it exists) that satisfies the particular case where:
\begin{equation} \label{particular-case}
\mbox{$q_j = q_k = q$, for any $1 \leq j,k\leq n_q$ and $n_p = n_q = n_r$.}
\end{equation}
 We get the following results about the factorization of $\sigma(T^2)$, for $T \in {\cal{D}}_A$.
\begin{theorem} \label{theoppal}
\emph{Let $A$ be a special perfect polynomial, denoted in (\ref{forme-de-A}), which satisfies (\ref{particular-case}).
%such that $q_j = q_k = q$, for any $1 \leq j,k \leq n_q$ and $n_p=n_q=n_r$.
Then: \\
i) For any $P \in {\cal{M}}_{in}$, $\sigma(P^2)$ is irreducible so that ${\cal{M}}_{ed} = {\cal{S}}_{1}$ and  $q=2p$.\\
ii) For any $Q \in {\cal{M}}_{ed}$, $\sigma(Q^2) = P \cdot R$, for a unique $(P,R) \in {\cal{M}}_{in} \times {\cal{M}}_{ax}$, so that $r = 3p$.\\
iii) For any $S \in {\cal{M}}_{ax}$, $\sigma({S}^2) = P \cdot Q \cdot R$, for a unique $(P,Q,R) \in {\cal{M}}_{in} \times {\cal{M}}_{ed} \times {\cal{M}}_{ax}$} or
$\sigma({S}^2) = P_1\cdots P_4 \cdot Q$ or $\sigma({S}^2) = {P_1}^2 P_2 P_3 \cdot Q$, for $P_j \in {\cal{M}}_{in}$ and $Q \in {\cal{M}}_{ed}$.
%iii) For any $S \in \varphi_1({\cal{M}}_{ax})$, $\sigma({S}^2) = R \cdot Q \cdot P$, for a unique $(P,Q,R) \in %{\cal{M}}_{in} \times {\cal{M}}_{ed} \times {\cal{M}}_{ax}$.\\
%iv) For any $S \in {\cal{M}}_{ax} \setminus \varphi_1({\cal{M}}_{ax})$, $\displaystyle{\sigma({S}^2) = \prod_{i=1}^4 %{P_i}}$, for some $P_i \in {\cal{M}}_{in}$.}
%In this case,
%either $(Q$ divides both $RP+1$ and $R+P+1)$ or $(Q$ divides both $RP+R+1$ and $R+P)$.
\end{theorem}

\section{Useful facts}
According to (\ref{forme-de-A}), we suppose that
%\begin{equation} \label{formedeA}
$$A = {P_1}^2 \cdots P_{n_p}^{2} \ {Q_1}^2 \cdots Q_{n_q}^2 \ {R_1}^2 \cdots R_{n_r}^2,$$
%\end{equation}
is an odd special perfect polynomial, where $1\leq i,j\leq n_p$, $1 \leq k, \ell \leq n_r$:
$$\deg(P_i) = \deg(P_{j})  < \deg(Q_1) \leq \cdots \leq \deg(Q_{n_q}) < \deg(R_k) = \deg(R_{\ell}).$$
\begin{lemma}\label{exposant-max}
Let $m \in \N$ and $S \in {\cal{D}}_A$ be such that $S^m \| \sigma(A)$, then $ m\leq 2$.
\end{lemma}

We recall the following five lemmas stated in \cite{Gall-Rahav-Bdx}.
\begin{lemma} \label{mini-medium-maxi}
One has, for any $i \not= j$, $\deg(P_i) = \deg(P_j) =p$, $\deg(R_i) = \deg(R_j)=r$. Moreover,
$n_p$ is even and $n_q \geq n_r\geq 3.$
\end{lemma}
%\begin{remark}
%{\rm{We may have $r = 2p$ or $p < \deg(Q_j) < 2p$, for some $j$.}}
%\end{remark}
\begin{definition} \label{min-med-max}
{\emph{We call the $P_i$'s: minimal (resp. the $Q_j$'s:  medium , the $R_k$'s:  maximal)  prime divisors of $A$.}}
\end{definition}
We refer to Notations in (\ref{lesPQR}) and in (\ref{lesPQR2}).
\begin{lemma} \label{alldegre-even}
The integers $p$, $r$ and all $q_j$ are even.
\end{lemma}
\begin{lemma} \label{samedegree}
Let $T_1, T_2$ be two divisors of $A$ of the same degree such that $T_1$ divides $\sigma({T_2}^2)$. Then, $T_2$ does not divide $\sigma({T_1}^2)$.
\end{lemma}
\begin{lemma} \label{maxdegree}
For any $R  \in {\cal{M}}_{ax}$, there exists a unique couple $(S,T)$ such that $S \in {\cal{M}}_{ax}$, $T \in {\cal{M}}_{ed}$ and $R$  divides both $\sigma(S^2)$ and $\sigma(T^2)$.
\end{lemma}
\begin{lemma} \label{alldegreeven}
If $A = {T_1}^2 \cdots {T_n}^2 \in \F_2[x]$ is perfect over $\F_2$, with each $T_i$ irreducible, then $T_i$ is of even degree and $T_i \equiv 1 \mod x^2+x+1$.
\end{lemma}
\begin{corollary} \label{noMersenne}
No Mersenne irreducible polynomial divides a special perfect polynomial.
\end{corollary}
\begin{proof}
A Mersenne polynomial is of the form $M_{ab} = x^a(x+1)^b + 1$, for some $a,b \in \N^*$. If it divides  a special perfect polynomial $A$, then $M_{ab} \equiv 1 \mod x^2+x+1$. So,
$x^a(x+1)^b \equiv 0 \mod  x^2+x+1$, which is impossible.
\end{proof}
\begin{lemma} \label{sigTdeux}
If  $S,T \in \F_2[x]$ where $T$ is irreducible and $\deg(T) \leq 2 \deg(S)$, then $S^2$ does not divide $\sigma(T^2)$.
%ii) For any $P \in {\cal{M}}_{in}$, $\sigma(P^2)$ is square-free.
\end{lemma}
\begin{proof}
If $S^2$ divides $\sigma(T^2) = 1+T+T^2$, then it also divides $T'$. So, $\deg(T) > \deg(T') \geq \deg(S^2) = 2 \deg(S)$.
%ii) follows from i) with $T = P$, since $\deg(P)$ is minimal.
\end{proof}
\begin{lemma} \label{omega-P-mini-2}
Let  ${\cal{T}} \subset  {\cal{D}}_A$ and $T_i, T_j,  T_{i1}, T_{i2},  T_{j1}, T_{j2} \in {\cal{T}}$ be such that they have the same degree and $\sigma({T_i}^2) = T_{i1} \cdot T_{i2}$, $\sigma({T_j}^2) = T_{j1} \cdot T_{j2}$.
Then,  $T_{i1}$, $T_{i2}$,  $T_{j1}$ and $T_{j2}$ are all distinct.
\end{lemma}
\begin{proof}
$T_{i1}$ and $T_{i2}$ are distinct because  $\sigma({T_i}^2)$ is not a square. If $T_{i1} =  T_{j1} = T$, then
$T$ divides  $\sigma({T_i}^2) +  \sigma({T_j}^2) = (T_i + T_j) \ (T_i + T_j+1)$, which is impossible by considering degrees.
\end{proof}
\section{The maps $\varphi_1, \varphi_2, \varphi_3$ and  $\varphi_4$} \label{lesphij}
From Lemma \ref{maxdegree}, we may define the following four maps:
$$\varphi_1 :  {\cal{M}}_{ax} \to {\cal{M}}_{ax}, \ \varphi_2 :  {\cal{M}}_{ax} \to {\cal{M}}_{ed}, \ \varphi_3 :  {\cal{M}}_{ax} \to \F_2[x], \  \varphi_4 :  {\cal{M}}_{ax} \to \F_2[x],$$
such that for any $R \in {\cal{M}}_{ax}$, one has
$$\sigma((\varphi_1(R))^2) = R \cdot \varphi_3(R),\ \sigma((\varphi_2(R))^2) = R \cdot \varphi_4(R).$$
Note that $\varphi_3(R)$ and $\varphi_4(R)$ may be reducible.

%For $f : {\cal{M}}_{ax} \to {\cal{M}}_{ax}$, we put:
%$$\mbox{$f^{(0)} = {\rm{Id}}_{{\cal{M}}_{ax}}$ and $f^{(m)} := f \circ \cdots \circ f$ ($m$-times), $m \geq 1$.}$$
\begin{lemma}
For any $R \in {\cal{M}}_{ax}$, $\varphi_4(R)$ factors in ${\cal{M}}_{in} \cup {\cal{M}}_{ed}.$
\end{lemma}
\begin{proof}
For $R \in {\cal{M}}_{ax}$,  $\sigma((\varphi_2(R))^2) = R \cdot \varphi_4(R)$. So,
$2 q \geq 2 \deg(\varphi_2(R)) = r + \deg(\varphi_4(R))$. Hence, $\deg(\varphi_4(R)) \leq 2q -r < r$.
\end{proof}
\begin{corollary} \label{phi2-inject}
The map $\varphi_2$ is injective.
\end{corollary}
\begin{proof}
%It suffices to prove that it is injective.
If  $\varphi_2(R_1) =  \varphi_2(R_2) = Q$, then $R_1 \cdot  \varphi_4(R_1) = \sigma(Q^2) = R_2 \cdot \varphi_4(R_2)$, with
$\deg(\varphi_4(R_j)) < r =\deg(R_1) = \deg(R_2)$. Therefore, $R_1 = R_2$.
\end{proof}
\begin{lemma} \label{sortdesRmax}
For any $R \in {\cal{M}}_{ax}$, $R^2$ divides $\sigma(A)$.
\end{lemma}
\begin{proof}
For any $R \in {\cal{M}}_{ax}$, $R$ divides both  $\sigma(({\varphi_1}(R))^2)$ and $\sigma(({\varphi_2}(R))^2)$.
\end{proof}
\begin{lemma} \label{lesphi2}
i) $\varphi_2(R) = \varphi_1(R) + R+1$, for $R \in {\cal{M}}_{ax}$\\
ii)  $\varphi_4(R) = \varphi_3(R)+R+1$,  for $R \in {\cal{M}}_{ax}$.
\end{lemma}
\begin{proof}
The polynomial $R$ divides both  $\sigma((\varphi_1(R))^2)$ and  $\sigma((\varphi_2(R))^2)$. Moreover, one has
$(\varphi_1(R) + \varphi_2(R)) \cdot  (\varphi_1(R) + \varphi_2(R)+1)=\sigma((\varphi_1(R))^2) + \sigma((\varphi_2(R))^2) = R \cdot (\varphi_3(R) + \varphi_4(R))$. Therefore, by considering degrees, $R$ being irreducible, one must have: $R=\varphi_1(R) + \varphi_2(R)+1$ and  $\varphi_3(R) + \varphi_4(R) =  \varphi_1(R) + \varphi_2(R) = R+1$.
\end{proof}

\section{Factorization of $\sigma(T^2)$ for $T \in {\cal{D}}_A$}
\subsection{Case where $T \in {\cal{M}}_{in}$} \label{factor-sigPmin-deux}
\begin{lemma} \label{sigPmindeux}
For any $T \in {\cal{M}}_{in}$, for any $S \in {\cal{D}}_{A}$,  $S^2$ does not divide $\sigma(T^2)$.
\end{lemma}
\begin{proof}
It follows from Lemma \ref{sigTdeux}, since $T$ is of minimal degree.
\end{proof}
\begin{lemma} \label{omega-P-mini}
If  $P \in {\cal{M}}_{in}$, then either $\sigma(P^2)$ is irreducible or $\sigma(P^2) = P_iP_j$ where $P_i, P_j \in {\cal{M}}_{in}$.
%ii) If $R \in {\cal{M}}_{ax}$ then $\sigma(R^2)$ is reducible. Moreover, if $\omega(\sigma(R^2)) = 2$  then $\sigma(R^2) = %R_1R_2$, where $R_1, R_2 \in {\cal{M}}_{ax}$.
\end{lemma}
\begin{proof}
One has $\omega(\sigma(P^2)) \leq 2$, because $\deg(P)$ is minimal. Our result follows by considering degrees.
%ii) Analogous argument: if $\sigma(R^2)= T$ is irreducible, then we get the contradiction:  $2r = \deg(T) \leq r$.
\end{proof}
\begin{comment}
\begin{remark} \label{factor-sigSdeux}
{\rm{If $T, T_1$ and $T_2$ are irreducible of the same degree and if they satisfy Lemma \ref{alldegreeven}, then $\sigma(T^2)$ may be irreducible or it may factor as $T_1 T_2.$\\
Example 1: One has  $\sigma(T^2) = T_1 \cdot T_2$, where
$$\begin{array}{l}
T = 1+x^3(x+1)^5(x^2 + x + 1)^3 ,\\
T_1 = 1+x(x+1)(x^2 + x + 1)(x^5 + x^4 + x^3 + x^2 + 1)(x^5 + x^3 + x^2 + x + 1),\\
T_2 = 1+x(x+1)(x^2 + x + 1)(x^{10} + x^9 + x^8 + x^6 + x^4 + x^3 + 1).
\end{array}$$
Exemple 2: $\sigma(T^2)$ is irreducible if $$T = 1+x(x+1)(x^2 + x + 1)(x^4 + x^3 + 1)(x^6 + x^5 + x^2 + x + 1).$$}}
\end{remark}
\end{comment}
\begin{lemma}
If $P \in  {\cal{M}}_{{in}1}$, then $P + \alpha$ and $P+\alpha+1$ are both irreducible over $\F_4$.
\end{lemma}
\begin{proof}
One has $\sigma(P^2) = 1+P+P^2 = (P+\alpha)(P+\alpha+1)$. If $P+\alpha$ is reducible, then $P+\alpha = S_1S_2$ where $0<\deg(S_1) < \deg(P)$.
Hence, $\sigma(P^2) =  (P+\alpha) (\overline{P+\alpha}) = S_1 \overline{S_1}  S_2 \overline{S_2}$, where $S_j \overline{S_j} \in \F_2[x]$. So, $\sigma(P^2)$ is reducible over $\F_2$. It contradicts the fact: $P \in  {\cal{M}}_{{in}1}$.
\end{proof}
\begin{lemma} \label{Q-divise-sigSdeux}
Let $P \in {\cal{M}}_{{in}1}$ and $S \in {\cal{D}}_A$ be such that $Q=\sigma(P^2)$ divides $\sigma(S^2)$.
Then, $\deg(S) > 2 \deg(P)$.
% and either $(Q$ divides both $SP+1$ and $S+P+1)$ or $(Q$ divides both $SP+S+1$ and $S+P)$.
\end{lemma}
\begin{proof}
First,  $Q=1+P+P^2 \in {\cal{S}}_{1}$ (with $P \in {\cal{M}}_{{in}1}$).
If $1+S+S^2 = (1+P+P^2) \cdot U$, then one has, over $\F_4$: $$(S+\alpha)(S+\alpha+1) = (P+\alpha) (P+\alpha+1) (V + \alpha W)(V+(\alpha+1)W),$$
where  $P+\alpha$ and  $P+\alpha+1$ are irreducible over $\F_4$,  $V,W \in \F_2[x],$ $\deg(V) > \deg(W) \geq 1$.
Therefore, w.l.o.g, we get:
$$S+\alpha = (P+\alpha) (V+ \alpha W) = PV + W + \alpha \cdot (V+PW+W),$$
or
$$S+\alpha = (P+\alpha+1) (V+ \alpha W) = PV +V + W + \alpha \cdot (V+PW),$$\\
Hence, $(PV + W = S$ and $V+(P+1)W = 1)$ or $((P+1)V+ W = S$ and $V+PW = 1)$.\\
By considering degrees, one has:
$$\deg(S) = \deg(P) + \deg(V),\ \deg(P) + \deg(W) = \deg(V).$$
Hence, $\deg(S) = 2 \deg(P) + \deg(W) > 2 \deg(P)$.
\begin{comment}
Now, if $(PV + W = S$ and $V+(P+1)W = 1)$, then $V = \displaystyle{\frac{PS + S + 1}{Q}}$ and $W = \displaystyle{\frac{S + P}{Q}}$. Therefore, $Q$ divides both $PS+S+1$ and $S+P$.\\
If $(P+1)V+ W = S$ and $V+PW = 1$, then $V = \displaystyle{\frac{PS+ 1}{Q}} $ and $W = \displaystyle{\frac{S + P+ 1}{Q}}$. Therefore, $Q$ divides both $SP+1$ and $S+P+1$.
\end{comment}
\end{proof}
\begin{lemma} \label{np-mp}
If ${\cal{M}}_{{in}2} \not= \emptyset$, then $\# {\cal{M}}_{{in}} \geq 4$.
\end{lemma}
\begin{proof}
Let $P \in {\cal{M}}_{{in}2}$. Then $\sigma({P}^2) = TV$, where $T, V \in {\cal{M}}_{{in}}$, $P \not=T \not= V$. Thus,
$n_p = \# {\cal{M}}_{{in}} \geq 4$, $n_p$ being even (Lemma \ref{mini-medium-maxi}).
\end{proof}
In the case where ${\cal{M}}_{{in}2} \not= \emptyset$, we may put $n_p = 2m_p \geq 4$ and we get:
\begin{lemma} \label{omega-P-mini-2-bis}  If ${\cal{M}}_{{in}2} \not= \emptyset$, then, after re-indexing:\\
i) ${\cal{M}}_{{in}2} = \{P_i \in {\cal{M}}_{in}: 1\leq i \leq m_p-1\}$, $\sigma({P_i}^2)  = P_{2i} \cdot P_{2i+1}$,\\
ii) for $m_p \leq i \leq 2m_p$, $\sigma({P_i}^2)$ is irreducible, so that $n_r+n_q \geq m_p + 1$.
\end{lemma}
\begin{proof}
i) - ii): By Lemma \ref{omega-P-mini-2}, with  ${\cal{T}} = {\cal{M}}_{in}$, for any $P,S,T \in {\cal{M}}_{in}$ pairewise distinct, $P$ does not divide simultaneously $\sigma(S^2)$ and $\sigma(T^2)$. Thus, after re-indexing, we get:
$$\sigma({P_1}^2) = P_2P_3,\ \sigma({P_2}^2) = P_4P_5, \ldots, \sigma({P_{m_p-1}}^2) = P_{2m_p-2}P_{2m_p-1}.$$
For $m_p \leq i \leq 2m_p$, $\sigma({P_i}^2) \in {\cal{M}}_{ed} \cup  {\cal{M}}_{ax}$. So, $m_p+1 \leq \#({\cal{M}}_{ed} \cup  {\cal{M}}_{ax}) = n_q + n_r$.
\end{proof}
\begin{corollary} \label{2p-inf-r} i)  ${\cal{M}}_{{in}1} \not= \emptyset$ so that ${\cal{S}}_{1} \not= \emptyset$.\\
ii) $2p < r < 2q$, \\
iii) ${\cal{S}}_{1} = \{\sigma(P^2): P \in {\cal{M}}_{{in}1}\} \subset {\cal{M}}_{ed}$.\\
iv) If $n_p=n_q$, then  ${\cal{S}}_{1} = {\cal{M}}_{ed}$ and $q = 2p$.
% and $\# {\cal{S}}_{1} =  m_p + 1 \geq 1$.
%If ${\cal{M}}_{{in}2} \not= \emptyset$, then one has: \\
\end{corollary}
\begin{proof}
i): If ${\cal{M}}_{{in}2} = \emptyset$, then ${\cal{M}}_{{in}1} = {\cal{M}}_{{in}}$. If ${\cal{M}}_{{in}2} \not= \emptyset$, then apply Lemma \ref{omega-P-mini-2-bis}-ii).\\
ii): First, $2p \leq r$ by definition of $r$. Let $P \in {\cal{M}}_{{in}1}$ and put $Q=1+P+P^2$ (which is irreducible and divides $\sigma(P^2)$). If $2p = r$, then for any $T \in  {\cal{D}}_A$, $\deg(T) \leq 2p$. Thus, by Lemma \ref{Q-divise-sigSdeux}, $Q$ does not divide $\sigma(T^2)$ if $T \not= P$. Hence, $Q \| \sigma(A)$. It contradicts the fact: $Q^2$ divides $A = \sigma(A)$.\\
Now, from Lemma \ref{maxdegree}, we immediately get: $r < 2q$.\\
iii): Since $2p < r$, $\sigma({P}^2) \not\in {\cal{M}}_{ax}$ for any $P \in {\cal{M}}_{{in}1}$. So, $\sigma({P}^2) \in {\cal{M}}_{ed}$.\\
iv) follows from iii) and from the fact: $\# {\cal{S}}_{1} = n_p = n_q = \# {\cal{M}}_{ed}$.
% and $\# {\cal{S}}_{1} = \# {\cal{M}}_{{in}1} = m_p + 1$.
\end{proof}
%In the rest of the paper, we suppose that $n_p = n_q = n_r =n$.
%we set ${\cal{M}}_{{ed}1}: = {\cal{S}}_{1} = \{\sigma(P^2): P \in  {\cal{M}}_{{in}1}\}$ and
%$$\mbox{${\cal{M}}_{{in}2} = \{P_i: 1 \leq i \leq m_p-1\}$, ${\cal{M}}_{{in}1} = \{P_i: m_p \leq i \leq 2m_p\}$,  if %${\cal{M}}_{{in}2} \not= \emptyset$.}$$

\subsection{Case where $T  \in {\cal{M}}_{ed}$} \label{factor-Q-med}
We get the following facts from Lemmas \ref{sigTdeux} and  \ref{Q-divise-sigSdeux}:
\begin{corollary} \label{Pdeux-ne-divise-pas}
If $T \in  {\cal{M}}_{ed}$ and  $P \in  {\cal{M}}_{in}$, then $P^2$ does not divide $\sigma(T^2)$.
\end{corollary}
\begin{corollary} \label{Qspec-ne-divise-pas}
If $P \in  {\cal{M}}_{{in}1}$ and $T \in  {\cal{M}}_{ed}$ such that $\deg(T) \leq 2p$, then
$Q=1+P+P^2$ does not divide $\sigma(T^2)$.
\end{corollary}
\begin{lemma}
Let $Q \in {\cal{M}}_{ed}$ and $R  \in {\cal{M}}_{ax}$. Then,
$R$ divides $\sigma(Q^2)$ if and only if $Q = \varphi_2(R)$. In this case, $R$ also divides $\sigma((\varphi_1(R))^2)$.
\end{lemma}
\begin{proof}
The sufficiency is immediate. For the necessity,
If $Q \not= \varphi_2(R)$, then $R$ divides the three distinct polynomials $\sigma(Q^2)$, $\sigma((\varphi_2(R))^2)$ and $\sigma((\varphi_1(R))^2)$. So, we get the contradiction~: $R^3$ divides $\sigma(A) = A$ and $R^2 \| A$.
\end{proof}
\subsection{Case where $T  \in {\cal{M}}_{ax}$} \label{factor-R-max}
We get from Lemmas \ref{sigTdeux} and \ref{Q-divise-sigSdeux}:
\begin{corollary}
i) If $Q \in {\cal{M}}_{ed}$ and $R \in {\cal{M}}_{ax}$, then $Q^2$ does not divide $\sigma(R^2)$.\\
ii)  If $R \in  {\cal{M}}_{ax}$ and  $P_1, P_2 \in  {\cal{M}}_{in}$, then ${P_1}^2  {P_2}^2$ does not divide $\sigma(R^2)$.
\end{corollary}
\begin{corollary} \label{Med1-divise-sigRdeux}
If ${\cal{M}}_{ed} = {\cal{M}}_{{ed}1}$, then for any $Q \in {\cal{M}}_{ed}$, there exists a unique $R \in {\cal{M}}_{ax}$ such that $Q$ divides $\sigma(R^2)$.
\end{corollary}
\begin{proof}
Set $Q = 1+P+P^2 = \sigma(P^2)$ which is an irreducible factor of $\sigma(P^2)$. Then, $Q$ divides $\sigma(S^2)$ for no
$S \in {\cal{M}}_{ed} \cup ({\cal{M}}_{in} \setminus \{P\})$. So, $Q$ must divide $\sigma(R^2)$, for some $R \in {\cal{M}}_{ax}$. \\
Now, we prove the uniqueness: if $Q$ divides both $\sigma({R_1}^2)$ and $\sigma({R_2}^2)$ with $R_1 \not= R_2$, then $Q^3$ divides $\sigma(A) = A$, which is impossible.
\end{proof}
\begin{lemma} \label{Q1Q2-sigRdeux}
If $Q_1, Q_2 \in {\cal{S}}_{1}$ and $R \in   {\cal{M}}_{{ax}}$, then the product $Q_1 Q_2$ does not divide $\sigma(R^2)$.
\end{lemma}
\begin{proof}
Put $Q_j=1+P_j+{P_j}^2 \in {\cal{S}}_{1}$, where $P_j \in {\cal{M}}_{{in}1}$.\\
If $1+R+R^2 = (1+P_1+{P_1}^2) \cdot (1+P_2+{P_2}^2) \cdot U$, then one has, over $\F_4$:
$$(R+\alpha)(R+\alpha+1) = (P_1+\alpha) (P_1+\alpha+1) (P_2+\alpha) (P_2+\alpha+1) (V + \alpha W)(V+(\alpha+1)W),$$
where  $P_j+\alpha$ and  $P_j+\alpha+1$ are irreducible over $\F_4$,  $V,W \in \F_2[x],$ $ \deg(P_1) = \deg(P_2) = \deg(V) > \deg(W) \geq 1$.
Therefore, w.l.o.g, we get four cases:
$$R+\alpha = \left\{ \begin{array}{l}
(P_1+\alpha)  (P_2+\alpha) (V+ \alpha W),\\
%= (P_1 P_2 + 1)V + (P_1+P_2+1)W + \alpha \cdot [(P_1+P_2+1)V+(P_1P_2+P_1+P_2)W],$$
(P_1+\alpha) (P_2+\alpha) (V+ (\alpha+1) W),\\
(P_1+\alpha) (P_2+\alpha+1) (V+ \alpha W),\\
(P_1+\alpha) (P_2+\alpha+1) (V+ (\alpha+1) W).
\end{array} \right.$$
Hence, in the first case, $R =  (P_1 P_2 + 1)V + (P_1+P_2+1)W$ and $1 = (P_1+P_2+1)V+(P_1P_2+P_1+P_2)W$.\\
By considering degrees, one has: $0 = \deg(1) = \deg(P_1P_2W)$, which is impossible.\\
We see by the same manner, that the other cases do not happen.
\end{proof}

\begin{remark}
\emph{It would be possible that $P^2$ divides $\sigma(R^2)$, for some $P  \in {\cal{M}}_{{in}}$ and  $R  \in {\cal{M}}_{{ax}}$.}
\end{remark}

\begin{comment}
\begin{corollary} \label{Med-one}
 One has: $\# {\cal{M}}_{{ed}1} =2m_p$ $($resp. $m_p+1)$ if ${\cal{M}}_{{in}2} = \emptyset$ $($resp. otherwise$)$.
\end{corollary}
\begin{proof}
If ${\cal{M}}_{{in}2} = \emptyset$, then $\# {\cal{M}}_{{ed}1} = \# {\cal{M}}_{{in}}$ and it equals  $\# {\cal{M}}_{{in}1}$, otherwise.
\end{proof}
\end{comment}

\begin{lemma}
Let $R, S  \in {\cal{M}}_{ax}$. Then,
$R$ divides $\sigma(S^2)$ if and only if $S = \varphi_1(R)$. In this case, $R$ also divides $\sigma((\varphi_2(R))^2)$.
\end{lemma}
\begin{proof}
The sufficiency is immediate. For the necessity,
If $S \not= \varphi_1(R)$, then $R$ divides $\sigma(S^2)$, $\sigma((\varphi_1(R))^2)$ and $\sigma((\varphi_2(R))^2)$. So, we get the contradiction~: $R^3$ divides $\sigma(A) = A$.
\end{proof}
\begin{lemma} \label{omega=2-et-phi1}
Let $R \in {\cal{M}}_{ax}$ be such that $\omega(\sigma(R^2)) = 2$. Then, $R \in  \varphi_1({\cal{M}}_{ax})$.
\end{lemma}
\begin{proof}
One has $\omega(\sigma(R^2)) = R_1 R_2$, where $R_1, R_2$ must lie in ${\cal{M}}_{ax}$.
We claim that $R = \varphi_1(R_1)$ or $R = \varphi_1(R_2)$. If $R \not= \varphi_1(R_1)$ and $R \not= \varphi_1(R_2)$, then, $R_1$ divides $\sigma((\varphi_1(R_1))^2)$, $\sigma((\varphi_2(R_1))^2)$ and $\sigma(R^2)$.
Thus, ${R_1}^3$ divides $\sigma(A) = A$, which is impossible.
\end{proof}

\section{The proof of Theorem \ref{theoppal}}
We refer to Notations in (\ref{forme-de-A}), (\ref{lesPQR}) and in (\ref{lesPQR2}).
We suppose that $A = \sigma(A)$, $q_j = q_k = q$ for any $1 \leq j,k\leq n_q$ and $n_p = n_q = n_r=n$.
Corollary \ref{2p-inf-r} gives
\begin{lemma} \label{sigP2irreduc}
For any $P \in {\cal{M}}_{{in}}$, $\sigma(P^2)$ is irreducible and
${\cal{M}}_{{ed}} = \{\sigma(P^2): P \in {\cal{M}}_{{in}}\}$, so that $q = 2p$.
\end{lemma}
\begin{definition} \label{def-lambda}
\emph{Corollary \ref{Med1-divise-sigRdeux} allows us to define the map
$$\mbox{$\lambda: {\cal{M}}_{{ed}} \to {\cal{M}}_{{ax}}$ such that $\lambda(Q) = R$ if and only if $Q$ divides $\sigma(R^2)$.}$$}
\end{definition}
\begin{lemma} \label{lambda-biject}
The map $\lambda$ is bijective,  i.e., for any $R \in  {\cal{M}}_{{ax}}$, there exists a unique $Q \in  {\cal{M}}_{{ed}}$ such that $Q \mid \sigma(R^2)$.
\end{lemma}
\begin{proof}
It suffices to prove the injectivity, since $n_q = n_r$. If $\lambda(Q_1) = \lambda(Q_2) = R$, with $Q_1 \not= Q_2$, then both $Q_1$ and $Q_2$ divide $\sigma(R^2)$. It contradicts Lemma \ref{Q1Q2-sigRdeux}.
\end{proof}

\begin{comment}
We suppose that $q=2p$ and $r=3p$. We consider the extremal case that is:
$$\mbox{$(n_r = n_q = \displaystyle{\frac{n_p}{2} + 1},$ if ${\cal{M}}_{{in}2} \not= \emptyset)$ and $(n_r=n_q = n_p,$ if ${\cal{M}}_{{in}2} = \emptyset)$.}$$
%We consider the following cases:
%$$\mbox{${\cal{M}}_{{in}2} \not= \emptyset$, ${\cal{M}}_{{in}2} = \emptyset$, $\varphi_3({\cal{M}}_{ax}) \subset  %{\cal{M}}_{ax}$...}$$
%In each case, we distinguish:
%$$\mbox{(I): $q_j =q$ for any $j$, (II): $q_j \not=q_k$ for some $j,k$.}$$
\subsection{Case ${\cal{M}}_{{in}2} \not= \emptyset$ and $n_r = n_q = \displaystyle{\frac{n_p}{2} + 1}$} \label{Min2-non-vide}
In this case, $\# {\cal{M}}_{{in}} \geq 4$, $n_p = 2m_p$, ${\cal{M}}_{{in}2} = \{P_j: 1 \leq j \leq m_p-1\}$ and
${\cal{M}}_{{in}1} = \{P_j: m_p \leq j \leq 2m_p\}$.\\
%\subsubsection{Case where $q_j = q$ for any $1 \leq j \leq n_q$} \label{Qsamedegree}
%By Lemma \ref{mini-medium-maxi} and Corollary \ref{2p-inf-r}, one has
%$r > 2p = q$ so that  $2q = 4p$.
%and
%$$\begin{array}{l}
%n_q \geq n_r \geq 3 \mbox{ and $n_q \geq m_p +1$ because $2p < r$}.
%\end{array}$$
Recall that $\sigma({P_i}^2) = Q_i \in {\cal{M}}_{{ed}1}$, for $m_p \leq i \leq 2m_p$.
\end{comment}
\begin{lemma}
One has:\\
i) $\varphi_4({\cal{M}}_{ax}) \subset {\cal{M}}_{in}$,
so that $r = 3p$,\\
ii) $\varphi_2$ is bijective, so that $\varphi_2({\cal{M}}_{ax}) = {\cal{M}}_{ed}$.
\end{lemma}
\begin{proof}
i): For any $R \in {\cal{M}}_{ax}$, $\sigma((\varphi_2(R))^2) = R \cdot \varphi_4(R)$. If $\varphi_4(R) \in {\cal{M}}_{ed}$ for some $R$, then by considering degrees, one has: $2q = r +q$ so that $q=r$, which is impossible. Thus, $\varphi_4({\cal{M}}_{ax}) \subset {\cal{M}}_{in}$ and
$4p=2q = r+p$.\\
ii): By Corollary \ref{phi2-inject}, $\varphi_2$ is injective. Thus, it is bijective, since $n_r = n_q$.
\end{proof}
\begin{comment}
\begin{remark} \label{phi4-inject-ou-non}
\emph{The map $\varphi_4$ may not be injective: it happens that $\varphi_4(R_1) = \varphi_4(R_2) = P$ with $R_1 \not= R_2$.
Indeed, putting $Q_j = \varphi_2(R_j)$, we may have: $R_1P = \sigma({Q_1}^2) \not= \sigma({Q_2}^2) = R_2P$.}
\end{remark}
\end{comment}
\begin{lemma}
%i)
For any $Q \in {\cal{M}}_{ed}$, there exists a unique $R \in {\cal{M}}_{ax}$ such that $Q$ divides $\sigma(R^2)$.
%ii) $\varphi_3({\cal{M}}_{ax}) \cap {\cal{M}}_{ax} = \emptyset$.
%iii) $\varphi_1({\cal{M}}_{ax}) = {\cal{M}}_{ax}$, so that $\varphi_1$ is bijective.
\end{lemma}
\begin{proof}
We recall that any $Q \in {\cal{M}}_{ed}$ is of the form $\sigma(P^2)$, with $P \in {\cal{M}}_{in}$.
So, our result follows from Corollary \ref{Med1-divise-sigRdeux}.
%$Q \| \sigma(P^2)$, $Q$ does not divide $\sigma(T^2)$ for any $T \in {\cal{M}}_{ed} \cup ({\cal{M}}_{in} \setminus %\{P\})$. Thus, $Q$ must divide $\sigma(R^2)$ for some $R \in {\cal{M}}_{ax}$. Unicity : voir avant...\\
%ii): If $\varphi_3(R) = S \in {\cal{M}}_{ax}$, then $\sigma((\varphi_1(R))^2) = R S$ is divisible by no $Q \in %{\cal{M}}_{ed}$, which contradicts i).
%iii): If there exists $S \in {\cal{M}}_{ax}$ such that $S \not= \varphi_1(R)$ for any $R \in {\cal{M}}_{ax}$, then %no $T \in {\cal{M}}_{ax}$ divides $\sigma(S^2)$. ???
\end{proof}

\begin{comment}
\begin{corollary} \label{div-precis}
For any $P_1, P_2 \in {\cal{M}}_{{in}}$, for any $Q_1, Q_2 \in {\cal{M}}_{{ed}}$, for any $R \in {\cal{M}}_{{ax}}$, ${P_1}^2$ does not divide $\sigma(R^2)$, ${Q_1}^2$ does not divide $\sigma(R^2)$, the product $Q_1 Q_2$ does not divide $\sigma(R^2)$.
\end{corollary}
\begin{proof}
It remains to prove that no ${T}^2$ divides $\sigma(R^2)$, if $T \in {\cal{M}}_{{in}}$. First, any such $T$ already divides some $\sigma(Q^2)$, $Q \in {\cal{M}}_{{ed}}$. If $P^2$ divides $\sigma(Q^2)$ or $\sigma(R^2)$ for some $P$, $Q$ and $R$, then $P^3$ divides $\sigma(A)$, which is impossible.
\end{proof}
\end{comment}

\begin{lemma} \label{factor-sigQ2}
If $Q \in {\cal{M}}_{ed}$, then \\
i) there exists a unique $(T,P,R) \in {\cal{M}}_{in} \times {\cal{M}}_{in} \times {\cal{M}}_{ax}$ such that $\sigma(Q^2) = R \cdot P$, $Q = \sigma(T^2)$ and $P \not= T$\\
ii) $\sigma(Q^2)$ is not divisible by $P_1 \cdot P_2$ with $P_1, P_2 \in {\cal{M}}_{in}$.
\end{lemma}
\begin{proof}
i) Any $Q \in {\cal{M}}_{ed}$ is of the form $\varphi_2(R)$, for a unique $R \in {\cal{M}}_{ax}$. Moreover, $Q = \sigma(T^2)$,  with $T \in {\cal{M}}_{in}$. Put then $P= \varphi_4(R) \in {\cal{M}}_{in}$. If $P=T$, then $T$ divides $\sigma(Q^2) = 1+T+T^4$, which is impossible.\\
ii) follows from i).
\end{proof}

\begin{comment}
\begin{corollary} \label{phi4-presk-inject}
The map $\varphi_4$ is injective from ${\cal{M}}_{ax} \setminus {\varphi_4}^{-1}(\{P_1, P_{n_p}\})$ into ${\cal{M}}_{in}$.
%, so that $n_r \leq n_p+2$.
\end{corollary}
\begin{proof} First, for no $j \in \{1,\ldots, n_p\}$, the polynomials $P_1$ and $P_{n_p}$ divide $\sigma({P_j}^2)$.
Let $R_1, R_2 \in {\cal{M}}_{ax} \setminus {\varphi_4}^{-1}(\{P_1, P_{n_p}\})$ be such that $\varphi_4(R_1) = \varphi_4(R_2)$.
Put $Q_1 =\varphi_2(R_1)$ and $Q_2 = \varphi_2(R_2)$. If $R_1 \not= R_2$, then $Q_1 \not= Q_2$ because $\varphi_2$ is injective (Lemma \ref{phi2-inject}). One has: $\sigma({Q_1}^2) = R_1 \cdot \varphi_4(R_1)$,
$\sigma({Q_2}^2) = R_2 \cdot \varphi_4(R_2)$, $\varphi_4(R_1) = \varphi_4(R_2) = P \not\in \{P_{1}, P_{n_p}\}$. Therefore, $P^2$ divides $\sigma(A)$. Moreover,
by Lemma \ref{omega-P-mini-2-bis}, $P$ also divides $\sigma({P_i}^2)$ for some $i \in \{1, \ldots, m_p - 1\}$. It follows that $P^3$ divides  $\sigma(A) = A$, which is impossible.
\end{proof}
\end{comment}
\begin{lemma} \label{phi3-injectif}~\\
i) For any $R \in {\cal{M}}_{{ax}}$, $\varphi_3(R) = Q \cdot P$ for some $(P,Q) \in {\cal{M}}_{{in}} \times {\cal{M}}_{{ed}}$.\\
ii) The map $\varphi_3$ is injective.
\end{lemma}
\begin{proof}~\\
i): if $\varphi_3(R) = S \in {\cal{M}}_{{ax}}$, then $S$ divides $\sigma((\varphi_1(R))^2), \sigma((\varphi_1(S))^2)$ and $\sigma((\varphi_2(S))^2)$. Thus, $S^3$ divides $\sigma(A) = A$, which is impossible.\\
One has: $\sigma((\varphi_1(R))^2) = R \cdot \varphi_3(R)$, with $\deg(\varphi_3(R)) = 3p$. Since the map $\lambda$ is bijective (Lemma \ref{lambda-biject}), there exists $Q \in {\cal{M}}_{{ed}}$ such that $Q \mid \varphi_3(R)$. So, we must have:  $\varphi_3(R) = Q \cdot P$, where $P \in   {\cal{M}}_{{in}}$.\\
ii): if $\varphi_3(R_1) = Q \cdot P = \varphi_3(R_2)$ with $R_1 \not= R_2$, then
%$Q_1 = \varphi_2(R_1) \not= \varphi_2(R_2)= Q_2$ and
$\sigma((\varphi_1(R_1))^2) = R_1 \cdot Q \cdot P$, $\sigma((\varphi_1(R_2))^2) = R_2 \cdot Q \cdot P$. Since $Q = \sigma(T^2)$ for some $T \in {\cal{M}}_{{in}}$, we get the contradiction : $Q^3$ divides $\sigma(A) = A$.
\end{proof}
We get from  Lemma \ref{sigTdeux}, Corollary \ref{Qspec-ne-divise-pas} and Lemma \ref{Q1Q2-sigRdeux}:
\begin{corollary} \label{pqr}
For any $P \in {\cal{M}}_{{in}}$, for any $S \in {\cal{M}}_{{in}} \setminus \{P\}$, for any $T \in {\cal{M}}_{ed}$, $Q=1+P+P^2$ does not divide
the product $\sigma(S^2) \sigma(T^2)$, but it divides $\sigma(R^2)$ for some $R \in {\cal{M}}_{ax}$ $($depending on $P).$
\end{corollary}

\begin{corollary}
For any $R \in {\cal{M}}_{ax}$, $\omega(\sigma(R^2)) \geq 3$.
\end{corollary}
\begin{proof}
If $\omega(\sigma(R^2)) = 2$, for some $R$, then $R = \varphi_1(S)$, where $S$ and $\varphi_3(S)$ belong to ${\cal{M}}_{ax}$. Thus,
$\# \{R \in {\cal{M}}_{ax}: Q \mid \sigma(R^2), \mbox{ for some $Q \in {\cal{M}}_{ed}$}\} < n$. Moreover,
for any $Q \in {\cal{M}}_{ed}$, one has $Q = \sigma(P^2) \| \sigma(P^2)$ and $Q \| \sigma(R^2)$, for some $P$ minimal and $R$ maximal. So, there exist $T \in {\cal{M}}_{ax}$, $Q_1, Q_2 \in {\cal{M}}_{ed}$ such that $Q_1 Q_2$ divides $\sigma(T^2)$. It contradicts~Lemma \ref{Q1Q2-sigRdeux}.
\end{proof}
We get our theorem from Lemma \ref{sigP2irreduc}, Lemma \ref{factor-sigQ2} and from
\begin{corollary}
For any $R \in {\cal{M}}_{ax}$, $\sigma(R^2)$ factors as follow:\\
i) if $R = \varphi_1(S)$, then $\sigma(R^2) = S \cdot Q \cdot P$, with $P \in {\cal{M}}_{in}$ and $Q \in {\cal{M}}_{ed}$,\\
ii)  if $R \not\in  \varphi_1({\cal{M}}_{ax})$, then $\sigma(R^2) =  Q \cdot P_1 P_2P_3P_4$ or $\sigma(R^2) =  Q \cdot {P_1}^2 P_2P_3$, with $P_j \in {\cal{M}}_{in}$ and $Q \in {\cal{M}}_{ed}$.
\end{corollary}

\def\thebibliography#1{\section*{\titrebibliographie}
\addcontentsline{toc}
{section}{\titrebibliographie}\list{[\arabic{enumi}]}{\settowidth
 \labelwidth{[
#1]}\leftmargin\labelwidth \advance\leftmargin\labelsep
\usecounter{enumi}}
\def\newblock{\hskip .11em plus .33em minus -.07em} \sloppy
\sfcode`\.=1000\relax}
\let\endthebibliography=\endlist

\def\biblio{\def\titrebibliographie{References}\thebibliography}
\let\endbiblio=\endthebibliography

\newbox\auteurbox
\newbox\titrebox
\newbox\titrelbox
\newbox\editeurbox
\newbox\anneebox
\newbox\anneelbox
\newbox\journalbox
\newbox\volumebox
\newbox\pagesbox
\newbox\diversbox
\newbox\collectionbox
%--------------------------------------------
\def\fabriquebox#1#2{\par\egroup
\setbox#1=\vbox\bgroup \leftskip=0pt \hsize=\maxdimen \noindent#2}
%--------------------------------------------
\def\bibref#1{\bibitem{#1}

%\mbox{}\ignorespaces

\setbox0=\vbox\bgroup}
%--------------------------------------------
\def\auteur{\fabriquebox\auteurbox\styleauteur}
\def\titre{\fabriquebox\titrebox\styletitre}
\def\titrelivre{\fabriquebox\titrelbox\styletitrelivre}
\def\editeur{\fabriquebox\editeurbox\styleediteur}

\def\journal{\fabriquebox\journalbox\stylejournal}

\def\volume{\fabriquebox\volumebox\stylevolume}
\def\collection{\fabriquebox\collectionbox\stylecollection}
%--------------------------------------------
{\catcode`\- =\active\gdef\annee{\fabriquebox\anneebox\catcode`\-
=\active\def -{\hbox{\rm
\string-\string-}}\styleannee\ignorespaces}}
%--------------------------------------------
{\catcode`\-
=\active\gdef\anneelivre{\fabriquebox\anneelbox\catcode`\-=
\active\def-{\hbox{\rm \string-\string-}}\styleanneelivre}}
%--------------------------------------------
{\catcode`\-=\active\gdef\pages{\fabriquebox\pagesbox\catcode`\-
=\active\def -{\hbox{\rm\string-\string-}}\stylepages}}
%--------------------------------------------
{\catcode`\-
=\active\gdef\divers{\fabriquebox\diversbox\catcode`\-=\active
\def-{\hbox{\rm\string-\string-}}\rm}}
%--------------------------------------------
\def\ajoutref#1{\setbox0=\vbox{\unvbox#1\global\setbox1=
\lastbox}\unhbox1 \unskip\unskip\unpenalty}
%--------------------------------------------
\newif\ifpreviousitem
\global\previousitemfalse
\def\separateur{\ifpreviousitem {,\ }\fi}
%--------------------------------------------
\def\voidallboxes
{\setbox0=\box\auteurbox \setbox0=\box\titrebox
\setbox0=\box\titrelbox \setbox0=\box\editeurbox
\setbox0=\box\anneebox \setbox0=\box\anneelbox
\setbox0=\box\journalbox \setbox0=\box\volumebox
\setbox0=\box\pagesbox \setbox0=\box\diversbox
\setbox0=\box\collectionbox \setbox0=\null}
%--------------------------------------------
\def\fabriquelivre
{\ifdim\ht\auteurbox>0pt
\ajoutref\auteurbox\global\previousitemtrue\fi
\ifdim\ht\titrelbox>0pt
\separateur\ajoutref\titrelbox\global\previousitemtrue\fi
\ifdim\ht\collectionbox>0pt
\separateur\ajoutref\collectionbox\global\previousitemtrue\fi
\ifdim\ht\editeurbox>0pt
\separateur\ajoutref\editeurbox\global\previousitemtrue\fi
\ifdim\ht\anneelbox>0pt \separateur \ajoutref\anneelbox
\fi\global\previousitemfalse}
%--------------------------------------------
\def\fabriquearticle
{\ifdim\ht\auteurbox>0pt        \ajoutref\auteurbox
\global\previousitemtrue\fi \ifdim\ht\titrebox>0pt
\separateur\ajoutref\titrebox\global\previousitemtrue\fi
\ifdim\ht\titrelbox>0pt \separateur{\rm in}\
\ajoutref\titrelbox\global \previousitemtrue\fi
\ifdim\ht\journalbox>0pt \separateur
\ajoutref\journalbox\global\previousitemtrue\fi
\ifdim\ht\volumebox>0pt \ \ajoutref\volumebox\fi
\ifdim\ht\anneebox>0pt  \ {\rm(}\ajoutref\anneebox \rm)\fi
\ifdim\ht\pagesbox>0pt
\separateur\ajoutref\pagesbox\fi\global\previousitemfalse}
%--------------------------------------------
\def\fabriquedivers
{\ifdim\ht\auteurbox>0pt
\ajoutref\auteurbox\global\previousitemtrue\fi
\ifdim\ht\diversbox>0pt \separateur\ajoutref\diversbox\fi}
%--------------------------------------------
\def\endbibref
{\egroup \ifdim\ht\journalbox>0pt \fabriquearticle
\else\ifdim\ht\editeurbox>0pt \fabriquelivre
\else\ifdim\ht\diversbox>0pt \fabriquedivers \fi\fi\fi
.\voidallboxes}
%--------------------------------------------

\let\styleauteur=\sc
\let\styletitre=\it
\let\styletitrelivre=\sl
\let\stylejournal=\rm
\let\stylevolume=\bf
\let\styleannee=\rm
\let\stylepages=\rm
\let\stylecollection=\rm
\let\styleediteur=\rm
\let\styleanneelivre=\rm

\begin{biblio}{99}

\begin{bibref}{Canaday}
\auteur{E. F. Canaday} \titre{The sum of the divisors of a
polynomial} \journal{Duke Math. J.} \volume{8} \pages 721-737 \annee
1941
\end{bibref}

\begin{bibref}{Gall-Rahav-F4}
\auteur{L. H. Gallardo, O. Rahavandrainy} \titre{On
perfect polynomials over $\F_4$}
\journal{Port. Math. (N.S.)} \volume{62(1)} \pages 109-122 \annee
2005
\end{bibref}

\begin{bibref}{Gall-Rahav-Bdx}
\auteur{L. H. Gallardo, O. Rahavandrainy} \titre{Odd
perfect polynomials over $\F_2$}
\journal{J. Th\'eor. Nombres Bordeaux} \volume{19} \pages 165-174 \annee
2007
\end{bibref}

\begin{bibref}{Gall-Rahav8}
\auteur{L. H. Gallardo, O. Rahavandrainy} \titre{All perfect
polynomials with up to four prime factors over $\F_4$}
\journal{Math. Commun.} \volume{14(1)} \pages 47-65 \annee 2009
\end{bibref}

\begin{bibref}{Gall-Rahav-up-Mersenne}
\auteur{L. H. Gallardo, O. Rahavandrainy} \titre{All even (unitary) perfect polynomials over $\F_2$ with only Mersenne primes as odd divisors} \journal{Kragujevac J. Math.}
\volume{49(4)} \pages 639-652 \annee 2025
\end{bibref}

\end{biblio}

\end{document}